\documentclass[oneside,notitlepage,12pt]{article}

\usepackage{amssymb}
\usepackage[leqno]{amsmath}
\usepackage{amsfonts}
\usepackage{amsopn}
\usepackage{amstext}
\usepackage{amsthm}

\usepackage[colorlinks, backref]{hyperref}

\providecommand{\cal}{\mathcal}
\renewcommand{\Bbb}{\mathbb}

\newenvironment{pf}{\begin{proof}}{\end{proof}}

\newcommand{\Aaa}{{\cal{A}}}
\newcommand{\Bee}{{\cal{B}}}

\newcommand{\Ef}{{\cal{F}}}

\newcommand{\Kay}{{\cal{K}}}

\newcommand{\Pee}{{\cal{P}}}

\newcommand{\Tau}{{\cal{T}}}
\newcommand{\Yu}{{\cal{U}}}
\newcommand{\Vee}{{\cal{V}}}
\newcommand{\Wu}{{\cal{W}}}

\newcommand{\Be}{{\Bbb{B}}}

\newcommand{\sig}{\sigma}

\renewcommand{\phi}{\varphi}
\renewcommand{\rho}{\varrho}

\newcommand{\ntr}{n\in\omega}

\newcommand{\loe}{\leqslant}
\newcommand{\goe}{\geqslant}

\newcommand{\subs}{\subseteq}
\newcommand{\sups}{\supseteq}
\newcommand{\nnempty}{\ne\emptyset}

\newcommand{\cl}{\operatorname{cl}}
\newcommand{\Int}{\operatorname{int}}
\newcommand{\w}{\operatorname{w}}

\newcommand{\id}{\operatorname{id}}

\newcommand{\dom}{\operatorname{dom}}
\newcommand{\rng}{\operatorname{rng}}

\newcommand{\proves}{\vdash}

\newcommand{\poset}{{\Bbb{P}}}

\newcommand{\meet}{\cdot}
\newcommand{\Meet}{\prod}
\newcommand{\join}{+}
\renewcommand{\Join}{\sum}
\newcommand{\Land}{\;\&\;}

\newcommand{\by}{/}

\newtheorem{tw}{Theorem}[section]
\newtheorem{wn}[tw]{Corollary}
\newtheorem{lm}[tw]{Lemma}
\newtheorem{prop}[tw]{Proposition}

\theoremstyle{definition}

\theoremstyle{remark}

\newcommand{\setof}[2]{\{#1\colon #2\}}

\newcommand{\seq}[1]{\langle #1 \rangle}

\newcommand{\sett}[2]{\{#1\}_{#2}}
\newcommand{\sn}[1]{\{#1\}} % singleton
\newcommand{\dn}[2]{\{#1,#2\}} % doubleton
\newcommand{\pair}[2]{\langle #1, #2 \rangle} % pair
\newcommand{\triple}[3]{\langle #1, #2, #3 \rangle} % triple
\newcommand{\map}[3]{#1\colon #2 \to #3} % A function
\newcommand{\img}[2]{#1[#2]} % image of a set
\newcommand{\inv}[2]{{#1}^{-1}[#2]} % preimage of a set

\newcommand{\power}[1]{\Pee(#1)}

\newcommand{\fin}[1]{[#1]^{<\omega}}

\newcommand{\ctbl}{\ensuremath{{\aleph_0}}}

\newcommand{\iso}{\cong}

\newcommand{\ord}{\operatorname{ord}}

\newcommand{\bL}{\mathbb L}
\newcommand{\bK}{\mathbb K}
\newcommand{\Pf}[1]{\operatorname{Pf}(#1)}
\newcommand{\ult}[1]{\operatorname{Ult}(#1)}
\newcommand{\ultn}[2]{\operatorname{Ult^{#1}}(#2)}
\newcommand{\ulte}{\operatorname{Ult}}
\newcommand{\closed}[1]{\operatorname{Closed}(#1)} % the lattice of closed sets
\newcommand{\closedn}[2]{\operatorname{Closed}^{#1}(#2)} % the lattice of closed sets
\newcommand{\V}{{\mathbb V}}
\newcommand{\W}{{\mathbb W}}

\newtheorem{problem}{Problem}

\title{Compact spaces, lattices, and absoluteness:\\ a survey}
\author{
{\sc Wies{\l}aw Kubi\'s}\\ \\
{\small Institute of Mathematics, Academy of Sciences of the Czech Republic}\\
\texttt{kubis@math.cas.cz}
}

\begin{document}
\maketitle

\begin{abstract}
Given a compact space in a fixed universe of set theory, one can naturally define its interpretation in any ZFC extension of the universe. We investigate the stability of some classes of compact spaces with respect to extensions of this sort. We show that the class of Eberlein/Gul'ko compacta is stable (= absolute). On the other hand, there are examples of Corson compacta which are no longer Corson in some forcing extensions.

All the material comes from the author's notes written between 2004 -- 2006.
\end{abstract}

\tableofcontents

\section{Preliminaries}

Assume $\pair X\Tau$ is a compact space in a universe of set theory $\V$ and assume that $\W$ is an extension of $\V$ satisfying the axioms of ZFC (or just a large enough part of ZFC). Clearly $X$ is again a Hausdorff topological space and $\Tau$ is its natural open base. Call this new topology $\Tau^\W$. In many cases, $\pair X{\Tau^\W}$ is no longer compact. One can ask whether there is a natural compact space which provides an interpretation of $X$ in the extension $\W$. This is indeed the case and one can see it for example starting with a compact $0$-dimensional space and looking at its Boolean algebra of clopen sets. Being a Boolean algebra is an absolute property and hence one can interpret the compact space as the space of all ultrafilters in a given universe. In the case of non-$0$-dimensional spaces it is natural to deal with lattices, although there are some technical details which make the situation less obvious.

One has to point out that stability of certain classes of compact spaces in forcing extensions had already been studied by Bandlow~\cite{Bandlow}, using another approach. Lattice-theoretic approach was essentially used by Okunev, Szeptycki and the author in \cite{KOSz}, inspired by a result of Todor\v cevi\'c \cite{TodorcevicBaireClassOne} saying that the class of Rosenthal compact spaces is absolute with respect to forcing extensions.
Our approach does not require the knowledge of forcing, we just deal with absoluteness between two standard (i.e. transitive) models of ZFC.

In Section \ref{kraty} we briefly describe the lattice-theoretic approach to compact spaces. Section \ref{rozszerzenia} deals with extensions of compact spaces and their basic properties. In Section \ref{Corson} we investigate absoluteness of some subclasses of the class of Corson compact spaces.

\subsection{Alexander's Subbase Lemma revisited}

In order to see a better motivation for the use of lattices in the theory of compact spaces, we present a ``point-less" version of Alexander's Lemma.
Recall that a filter $F$ is \emph{generated by} $S$ if $S \subs F$ and for every $a \in F$ there is a finite set $s \subs S$ such that $\Meet s \loe a$.
We use the symbols $\join$ and $\meet$ for the join (supremum) and the meet (infimum) in a lattice.

\begin{lm}\label{LemAlexptls}
Let $\bL$ be a lattice generated by $G$.
Then every ultrafilter $p \subs \bL$ is generated by $p \cap G$.
\end{lm}

\begin{pf}
Fix $a \in p \in \ult \bL$.
Then $a = s_0 \join \dots \join s_{n-1}$, where each $s_i$ is a finite meet of some elements of $G$.
Since $p$ is a maximal filter, it follows that $s_j \in p$ for some $j < n$.
\end{pf}

\begin{wn}[Alexander's Subbase Lemma]
Let $\Yu$ be an open subbase for a topological space $\pair X \Tau$ such that every cover consisting of sets in $\Yu$ has a finite subcover.
Then $\pair X \Tau$ is a compact space.
\end{wn}

\begin{pf}
Let $\Bee = \setof{X \setminus U}{U \in \Yu}$.
Then $\Bee$ is a closed subbase with the property that every centered subcollection of $\Bee$ has a nonempty intersection.
Let $\bL$ be the lattice of sets generated by $\Bee$.
Then $\bL$ is a closed base for the topology $\Tau$ and by Lemma~\ref{LemAlexptls} every centered subcollection of $\bL$ has a nonempty intersection.
This obviously implies compactness, because every centered collection of closed sets can be replaced by a subcollection of $\bL$, with the same intersection.
\end{pf}

\section{Compact spaces and lattices}\label{kraty}

In this section we review the relationship between compact spaces and normal distributive lattices. 

Given a compact Hausdorff space $X$, fix a lattice of closed sets $\bL$ which is at the same time a closed subbase for $X$. A natural question is whether $X$ can be recovered from $\bL$. On the other hand, given an abstract lattice $\bL$, one can ask when there exists a compact Hausdorff space $X$ such that $\bL$ is isomorphic to a generating sublattice of closed subsets of $X$. This section is devoted to these questions. We study normal lattices, for which the natural Wallman topology on the ultrafilters is Hausdorff.
We also show how homomorphisms translate to continuous maps and vice versa. 
This algebraic approach to compact spaces is quite useful when studying compact spaces in forcing extensions.

The results presented here are mostly a folklore, although it is difficult to find them written explicitly.

There also exists a duality, involving prime filters, between distributive lattices and certain compact $T_0$ spaces, although we omit this topic here (except in Theorem~\ref{v-beer} below). We also omit the presentation of Stone's duality between Boolean algebras and compact 0-dimensional spaces, since it is well known and follows from the more general theory of normal lattices.

\subsection{Normal separative lattices}
We start with basic definitions.

Recall that a {\em lattice} is a partially ordered set $\pair L\loe$ such that both $a\meet b:=\inf\dn ab$ and $a\join b:=\sup\dn ab$ exist for every $a,b\in L$ and there are elements $0,1\in L$ such that $0\loe x\loe 1$ for every $x\in L$ \footnote{The last condition is usually not required, our definition refers to {\em bounded lattices}.}. 
We shall consider lattices as algebraic objects and therefore we shall write $\bL=\seq{L,\join,\meet}$ instead of $\bL=\pair L\loe$.
A lattice $\bL$ is {\em distributive} if $a\meet(b\join c)=(a\meet b)\join (a\meet c)$ holds for every $a,b,c\in L$. Every distributive lattice is isomorphic to a lattice of sets, therefore we shall use notions specific to sets, e.g. elements $a,b$ are {\em disjoint} if $a\meet b=0$. 
A lattice $\bL$ is {\em normal} 
if it is distributive\footnote{We assume distributivity here for the sake of convenience only.} and for every disjoint elements $a,b\in L$ there are $a_1,b_1\in L$ such that $a_1\meet b=0=a\meet b_1$ and $a_1\join b_1=1$. A lattice $\bL$ is {\em separative}
if it has this property as a poset, i.e. if $0<a\not\loe b$ then there is $c>0$ such that $c\loe a$ and $c\meet b=0$.

Given a lattice $\bL=\seq{L,\join,\meet}$, the structure $\seq{L,\meet,\join}$ is again a lattice, which will be called the {\em opposite lattice}
of $\bL$.

A {\em filter}
in a lattice $\bL=\seq{L,\join,\meet}$ is a subset $F$ of $L$ such that $0\notin F$, $1\in F$ and $a,b\in F$ implies $\setof{x\in L}{x\goe a\meet b}\subs F$. A filter $F$ is {\em prime}
if $a\join b\in F$ implies $a\in F$ or $b\in F$. 
A filter in the opposite of $\bL$ is called an {\em ideal}
of $\bL$.
Observe that a filter is prime iff its complement is an ideal. We denote by $\Pf \bL$ the collection of all prime filters in $\bL$.
Let $\ult\bL$ denote the collection of all maximal filters
(i.e. {\em ultrafilters}) in $\bL$. It is easy to see that $\ult\bL\subs\Pf\bL$.
Note that if a filter $p\subs\bL$ is maximal then for every $a\notin p$ there exists $b\in p$ such that $a\meet b=0$. A subset $M$ of a lattice $\bL$ is {\em centered}
if $m_0\meet\dots\meet m_{n-1}>0$ for every $m_0,\dots,m_{n-1}\in M$. In that case the set $$[M):=\setof{x\in\bL}{(\exists\;m_0,\dots,m_{k-1}\in M)\;m_0\meet\dots\meet m_{k-1}\loe x}$$ is a filter (the filter {\em generated by} $M$).
An intersection of sublattices is a sublattice, therefore every subset of a lattice {\em generates} some lattice. We say that $M\subs \bL$ is {\em generating}
if $\bL$ is the lattice generated by $M$ or, in other words, $\bL$ is the only sublattice of $M$ which contains $M$. In case where $\bL$ is distributive, the sublattice generated by $M$ consists of all elements of the form $\Join_{i<m}\Meet_{j<n}a_{i,j}$, where $\setof{a_{i,j}}{i<m,\;j<n}\subs M$ (symbols $\Join$ and $\Meet$ denote the supremum and the infimum respectively).

A map of lattices $\map f\bK\bL$ is a {\em lattice homomorphism} 
if $f$ preserves both meet and join and $f(0)=0$, $f(1)=1$.

For the general lattice theory we refer to \cite{Birkhoff} and \cite{Gratzer}.

\begin{prop}\label{kbeer} Let $\bL$ be a normal lattice. Then every prime filter extends uniquely to an ultrafilter. More precisely, if $F$ is a prime filter then
$$p_F=\setof{a\in\bL}{(\forall\;b\in F)\; a\meet b>0}$$ is an ultrafilter. \end{prop}

\begin{pf} Let $F$ be a prime filter in $\bL$. Assume $p,q$ are ultrafilters such that $F\subs p\cap q$. Suppose $p\ne q$. There are $a,b$ such that $a\in p\setminus q$, $b\in q\setminus p$ and $a\meet b=0$. Using normality, find $a',b'$ such that $a\meet b'=0=a'\meet b$ and $a'\join b'=1$. Then either $a'\in F$ or $b'\in F$ and in both cases we get a contradiction.

To show the ``more precisely" part, denote by $p$ the unique ultrafilter extending $F$. If $a\in p$ then $a\meet b>0$ for every $b\in F$. If $a\notin p$ then $F\cup \sn a$ is not centered (otherwise it would extend to an ultrafilter) and hence $a\meet b=0$ for some $b\in F$.
\end{pf}

A classical result of Stone \& Birkhoff says that a lattice $\bL$ is distributive if and only if for every filter $F\subs\bL$ and for every $a\in \bL\setminus F$ there exists a prime filter $p$ such that $a\notin p$ and $F\subs p$. In order to make this section self-contained, we prove the crucial part of this result below.

\begin{prop} Let $\bL=\seq{L,\join,\meet}$ be a distributive lattice and let $I\subs L$ be an ideal. Then every maximal filter disjoint from $I$ is prime.
\end{prop}

\begin{pf} Let $F$ be a maximal filter disjoint from $I$. Extend $I$ to a maximal ideal $J$ disjoint from $F$. We need to show that $J\cup F=L$. Suppose this is not the case and fix $x\in L\setminus(J\cup F)$. Then the ideal generated by $J\cup\sn x$ intersects $F$, which means that $x\join j=:f_1\in F$ for some $j\in J$. By the same reason, $x\meet f=:j_1\in I$ for some $f\in F$. On the other hand we have
$$f_1\meet f = (x\join j)\meet f = (x\meet f)\join (j\meet f)=j_1\join (j\meet f)\loe j_1\join j.$$
Thus $F\cap J\nnempty$, a contradiction.
\end{pf}

Fix a lattice $\bL$. Given $a\in \bL$ define
$$a^-=\setof{p\in\Pf\bL}{a\notin p},\qquad a^+=\setof{p\in\Pf \bL}{a\in p}.$$
Note that $(a\meet b)^+=a^+\cap b^+$, $(a\join b)^+=a^+\cup b^+$ and $0^+=\emptyset$, $1^+=\Pf\bL$. Moreover $a^-=\Pf\bL\setminus a^+$. Thus the map $a\mapsto a^+$ is a lattice homomorphism from $\bL$ to $\seq{\power\bL,\cup,\cap}$ and it is one-to-one if $\bL$ is separative.
The set $\Pf\bL$ is naturally equipped with a topology which is generated by sets of the form $a^-$, where $a\in\bL$. Denote this topology by $\Tau_\bL$. Then $\setof{a^-}{a\in\bL}$ is an open base for $\Pf\bL$ which is at the same time a lattice isomorphic to the opposite of $\bL$. The subspace topology on $\ult\bL$ (which we shall also denote by $\Tau_\bL$) turns out to be compact too, even though in many cases $\ult\bL$ is dense in $\Pf\bL$.
The following result is rather known, although we were unable to find explicit reference.

\begin{tw}\label{v-beer} Let\/ $\bL$ be a distributive lattice. Then:
\begin{enumerate}
	\item[$(a)$] The space $\pair{\Pf\bL}{\Tau_\bL}$ is compact and $T_0$. It is not\/ $T_1$ unless $\bL$ is a Boolean algebra (in that case it is Hausdorff).
	\item[$(b)$] A closed subset of $\pair{\Pf\bL}{\Tau_\bL}$ equals $a^+$ for some $a\in\bL$ if and only if its complement is compact.
	\item[$(c)$] The space $\pair{\ult\bL}{\Tau_\bL}$ is compact and $T_1$. It is Hausdorff if and only if\/ $\bL$ is normal.
	\item[$(d)$] $\ult\bL=\Pf\bL$ if and only if\/ $\bL$ is a Boolean algebra.	
\end{enumerate}
\end{tw}

\begin{pf} In order to show compactness, fix $\Aaa\subs\bL$ such that the family $\setof{a^+}{a\in\Aaa}$ is centered in $\Pf\bL$ or $\ult\bL$. This is equivalent to the fact that $a_0\meet\dots\meet a_{k-1}>0$ for every $a_0,\dots,a_{k-1}\in\Aaa$. In particular, ``centered in $\Pf\bL$" implies ``centered in $\ult\bL$".
Define
$$F=\setof{x\in\bL}{(\exists\;a_0,\dots,a_{k-1}\in\Aaa)\; a_0\meet\dots\meet a_{k-1}\loe x}.$$
Then $F$ is a filter in $\bL$ and therefore there exists $p\in\ult\bL$ which extends $F$. We have $p\in\ult\bL\cap\bigcap_{a\in\Aaa}a^+$. This shows the compactness of both $\Pf\bL$ and $\ult\bL$. A similar argument shows that sets of the form $b^-$ are compact for every $b\in\bL$ (using the fact that an element is separated from a filter by a prime filter). Now assume $F\subs\Pf\bL$ is closed and $\Pf\bL\setminus F$ is compact. Using compactness and the fact that $\setof{a^-}{a\in\bL}$ is an open base, we get $\Pf\bL\setminus F=a_0^-\cup\dots\cup a_{k-1}^-$ for some $a_0,\dots,a_{k-1}\in\bL$. Thus $F=\bigcap_{i<k}a_i^+=(a_0\meet\dots\meet a_{k-1})^+$.
This shows (b).

For the proof of (c), fix $p\ne q$ in $\ult\bL$. Fix $a\in p\setminus q$ and using the maximality of $q$ find $b\in q\setminus p$ such that $a\meet b=0$. Thus $p\in a^+\setminus b^+$ and $q\in b^+\setminus a^+$, which shows that $\ult\bL$ is $T_1$.
If additionally $\bL$ is normal, the disjoint sets $a^+$, $b^+$ are separated by disjoint basic open sets and therefore $\ult\bL$ is Hausdorff.
Conversely, assume $\ult\bL$ is Hausdorff and fix disjoint $a,b\in\bL$. Then $a^+$ and $b^+$ are disjoint closed sets in $\ult\bL$. Using compactness, $T_2$ and the fact that $\setof{x^+}{x\in\bL}$ is a closed base, we can find $a_1,b_1\in\bL$ with $a^+\cap b_1^+=\emptyset=a_1^+\cap b^+$ and $a_1^+\cup b_1^+=\ult\bL$. Translating it back to the lattice, we get normality. This shows (c).

We now show (d). Assume every prime filter is maximal in $\bL$. We shall show that $\Pf\bL$ is Hausdorff, which by (b) will give us that each set of the form $a^-$ is closed and therefore, again by (b), it equals $b^+$ for some $b\in \bL$; necessarily $b$ is the complement of $a$ in $\bL$.

Fix $p\ne q$ in $\Pf\bL$. Fix $a\in p\setminus q$. 
Then $I=\bL\setminus p$ is a prime ideal which does not contain $a$. Let $J\sups I$ be a maximal ideal. Then $J$ is prime and $\bL\setminus J$ is a prime filter which is a subset of $p$. Thus $J=I$, i.e. $I$ is maximal. Hence $a\join b=1$ for some $b\in I$. Then $b\in q$, because $q$ is prime. Thus $a^-$, $b^-$ are disjoint neighborhoods of $q$ and $p$ respectively. This shows that $\Pf\bL$ is Hausdorff and finishes the proof of (d).

It remains to complete the proof of (a).
It is clear that $\Pf\bL$ is $T_0$. Assume it is $T_1$. Then $\ult\bL=\Pf\bL$. Indeed, if $p$ is a prime filter which is not maximal then there is $q\in\ult\bL$ with $p\subs q$ and $p\ne q$. In this case $p\in a^+$ implies $q\in a^+$ for every $a\in\bL$, i.e. no basic closed set separates $p$ from $q$, a contradiction. By (d), $\bL$ is a Boolean algebra.
\end{pf}

From now on, we shall consider spaces of ultrafilters only and we shall write $a^+$ and $a^-$ instead of $a^+\cap\ult\bL$ and $a^-\cap \ult\bL$.

\begin{wn} Every finite separative distributive lattice is a Boolean algebra. \end{wn}

\begin{pf} Let $\bL$ be a finite separative distributive lattice. By Theorem \ref{v-beer}(c), $\ult\bL$ is a finite $T_1$ space and hence it is a discrete space.
By separativity, $\bL$ is isomorphic to $\bL'=\setof{a^+}{a\in\bL}$ and $\bL'$ is a closed base for the topology of $\ult\bL$. Thus $\bL'=\power{\ult\bL}$.
\end{pf}

Let us note that there exist normal lattices which are far from being separative. For example, given a lattice $\bL$ we can add a new element $s$ and declare it to be strictly below every element of $\bL$; this new lattice, which can be denoted by $1+\bL$, is distributive and normal, because there are no nontrivial disjoint pairs of elements; on the other hand there is only one ultrafilter on $1+\bL$.

\subsection{Homomorphisms and continuous maps}

Let $\bK, \bL$ be normal lattices and let $\map h{\bK}{\bL}$ be a lattice homomorphism. One can expect that there exists a continuous map $\map f{\ult\bL}{\ult\bK}$ which is induced by $h$, i.e. $\inv f{a^+} = h(a)^+$ holds for every $a\in \bK$. This is not quite true in general, however we shall show that $h$ indeed induces a continuous map. 
If $p\in\ult\bL$ then $\inv hp$ is a prime filter in $\bK$ and therefore it extends uniquely to an ultrafilter, by Proposition \ref{kbeer}. This suggests the definition of a map $\map f{\ult\bL}{\ult\bK}$. We shall call $f$ the {\em map induced by} $h$ and we shall write $f=\ult h$ ($\ulte$ is indeed a functor on the category of normal lattices).

Recall that given a lattice $\bL$, we denote by $a^+$ the set of all {\em ultrafilters} in $\bL$ which contain $a$ (and $a^-$ is the complement of $a^+$ in $\ult\bL$).

\begin{lm}\label{lbeer} Let $\bK,\bL$ be normal lattices and let $\map h{\bK}{\bL}$ be a homomorphism. Then the map $\map f{\ult\bL}{\ult\bK}$ induced by $h$ (i.e. $f(p)$ is defined to be the unique ultrafilter extending $\inv hp$) is continuous.
\end{lm}

\begin{pf} Fix $p\in \ult\bL$ and $a\in\bK$ such that $f(p)\in a^-$. Then $\inv hp\cup \sn a$ is not centered, therefore there exists $b\in\bK$ such that $h(b)\in p$ and $a\meet b=0$. By normality, find $c,d\in \bK$ such that $a\meet d=0=b\meet c$ and $c\join d=1$. Then $h(b)\meet h(c)=0$ and hence $h(c)\notin p$. Thus $h(c)^-$ is a neighborhood of $p$. If $q\in h(c)^-$ then $h(d)\in q$, because $h(c)\join h(d)=1$; Hence $a\notin f(q)$, because $a\meet d=0$. It follows that $\img f{h(c)^-}\subs a^-$.
\end{pf}

We now describe homomorphisms which can be recovered from their induced continuous maps. Call a homomorphism $\map h\bK\bL$ {\em separative}
if for every $a\in \bK$, $b\in\bL$, $$h(a)\meet b=0\implies (\exists\;c\in\bK)\; a\meet c=0\Land b\loe h(c).$$

\begin{prop}\label{mbeer} Let $\map h\bK\bL$ be a homomorphism of normal separative lattices. The following conditions are equivalent:
\begin{enumerate}
	\item[$(a)$] $h$ is separative.
	\item[$(b)$] $\inv hp\in\ult\bK$ whenever $p\in\ult\bL$.
	\item[$(c)$] $\inv{\ult h}{a^+}=h(a)^+$ for every $a\in\bK$.
\end{enumerate}\end{prop}

\begin{pf} $(a)\implies(b)$ Assume $a\notin \inv hp$. Then $h(a)\meet b=0$ for some $b\in p$. By $(a)$, there is $c\in\bK$ such that $a\meet c=0$ and $b\loe h(c)$. In particular $c\in \inv hp$, so $\inv hp\cup\sn a$ is not centered. Hence $\inv hp$ is a maximal filter.

$(b)\implies(c)$ Let $f=\ult h$. We have
$$\inv f{a^+}=\setof{p\in\ult\bL}{a\in f(p)}=\setof{p\in\ult\bL}{h(a)\in p}=h(a)^+.$$

$(c)\implies(a)$ Assume $h(a)\meet b=0$. Then $b^+$ is a closed subset of $\ult\bL$ and therefore $\img f{b^+}$ is closed and disjoint from $a^+$, where $f=\ult h$. 
Using the fact that $\setof{x^+}{x\in\bK}$ is a closed base of $\ult\bK$, an easy compactness argument (see the proof of Lemma \ref{wbeer} below for details) gives an element $c\in\bK$ such that $a^+\cap c^+=\emptyset$ and $\img f{b^+}\subs c^+$. Then $a\meet c=0$ and $b^+\subs \inv f{c^+}$. By the separativity of $\bL$, we have $b\loe h(c)$.
\end{pf}

Let us note that if $\bK$ is a normal lattice which is not a Boolean algebra, then there exists a homomorphism $h$ on $\bK$ which is not separative. Indeed, by Theorem \ref{v-beer}(d), there exists a prime filter $F\subs\bK$ which is not maximal.
Let $\map h\bK{\bK\by F}$ be the quotient homomorphism. Then $\bK\by F\iso 2$ (the 2-element lattice) and $f=\ult h$ is an embedding of a singleton into $\ult\bK$. Let $q$ be the unique ultrafilter extending $F$ and choose $a\in q\setminus F$. Then $h(a)=0= h(a)\meet 1$ and there is no $c$ such that $a\meet c=0$ and $h(c)=1$, because $a\meet c=0$ iff $c\notin F$. It is not hard to see that a homomorphism induced by a continuous map of compact spaces is separative.

Let $\bL$ be a lattice and let $K\subs\bL$. We say that $K$ {\em separates} $\bL$
if for every disjoint $a,b\in\bL$ there exists $c\in K$ such that $a\loe c$ and $c\meet b=0$. In this case, there exists also $d\in K$ such that $b\loe d$ and $c\meet d=0$, i.e. every pair of disjoint elements of $\bL$ is separated by disjoint elements of $K$.

Summarizing the relationship between lattices and compact spaces, we get the following:

\begin{tw}\label{monoepik} $\ulte$ is a contravariant functor on the category of normal lattices into the category of compact Hausdorff spaces. Moreover, if $\map h\bK\bL$ is a homomorphism of normal lattices then
\begin{enumerate}
	\item[$(a)$] $h^{-1}(0)=\sn 0$ implies that $\ult h$ is a surjection;
	\item[$(b)$] $\img h\bK$ separates $\bL$ implies that $\ult h$ is an embedding.
\end{enumerate}
\end{tw}

\begin{pf} It is clear that $\ult{\id_\bL}=\id_{\ult\bL}$. We need to show that $\ult{g\circ h} = \ult h \circ \ult g$ for every $g,h$ such that $\rng h\subs\dom g$.
Fix $p\in\dom(\ult h)$. Let $\phi=\ult{g\circ h}$, $\rho=\ult g$, $\eta=\ult h$.
Then $\phi(p)$ is the unique ultrafilter extending $\inv{(g\circ h)}p$. We also have $\inv{(g\circ h)}p\subs \inv h{\rho(p)}\subs \eta(\rho(p))$. Thus $\phi(p)=\eta(\rho(p))$.

Now assume $\map h\bK\bL$ is a homomorphism such that $h^{-1}(0)=\sn0$.
Fix $p\in\ult\bK$. Then $\img hp$ is centered and hence extends to an ultrafilter $q\in\ult\bL$. Thus we have $p\subs\inv hq$ and therefore $p=\inv hq=\ult h(q)$.

Finally, assume $h$ is an epimorphism. Fix $p\ne q$ in $\ult\bL$. Then there are disjoint $a,b$ such that $a\in p$, $b\in q$. Since $\img h\bK$ separates $\bL$, there are $a',b'$ such that $a\loe h(a')$, $b\loe h(b')$ and $a'\meet b'=0$. Thus $a'\in\inv hp\setminus\ult h(q)$ and $b'\in\inv hq\setminus\ult h(p)$, which shows that $\ult h(p)\ne \ult h(q)$.
\end{pf}

%[How about the converse implications???]

\subsection{Compact spaces and generating lattices}

Let $X$ be a topological space. We denote by $\closed X$ the lattice of all closed subsets of $X$. Observe that $\closed X$ is a normal lattice iff $X$ is a normal topological space. Moreover, $\closed X$ is separative if $X$ is $T_1$. One can consider various sublattices of $\closed X$ which at the same time form a closed subbase. In case of a compact Hausdorff space $X$, any sublattice of $\closed X$ which is a closed subbase is normal and $X$ is homeomorphic to its space of ultrafilters.
We prove a more precise statement below.

We need one more notion, in order to simplify some statements: given a compact space $X$, a sublattice of $\closed X$ will be called {\em basic}
if it is a closed base for $X$. 

\begin{lm}\label{wbeer} Let $X$ be a compact space and assume $\bL\subs \closed X$ is a basic lattice. Then for every closed set $a$ and its neighbourhood $v$ there exists $b\in \bL$ such that $a\subs \Int b\subs b\subs v$.
\end{lm}

\begin{pf} By normality we can find a closed set $a'\subs v$ such that $a\subs \Int a'$. For each $x\in K\setminus v$ find $b_x\in\bL$ such that $K\setminus b_x$ is a neighborhood of $x$ disjoint from $a'$. By compactness, find $x_0,\dots, x_{k-1}\in K\setminus v$ such that $K\setminus v$ is covered by $(K\setminus b_{x_0})\cup\dots\cup (K\setminus b_{x_{k-1}})$. Let $b=b_{x_0}\cap\dots\cap b_{x_{k-1}}$. Then $a'\subs b\subs v$. \end{pf}

\begin{prop} Let $X$ be a Hausdorff compact space and let $\bL\subs\closed X$ be a sublattice which is a closed subbase for $X$. Then $\bL$ is normal, separative and $X$ is homeomorphic to $\ult\bL$.
\end{prop}

\begin{pf} Assume $a,b\in\bL$ are such that $a\not\subs b$ and pick $x\in a\setminus b$. By Lemma \ref{wbeer} there is $c\in \bL$ such that $x\in c$ and $c\cap b=\emptyset$. This shows separativity. A similar argument shows normality (using the fact that Hausdorff compact implies normal). Now, the embedding $\map j\bL{\closed X}$ is a lattice homomorphism, so it induces a map $f=\ult j$. Lemma \ref{wbeer} says that $j$ is a separative homomorphism. Clearly, $\ult{\closed X}=X$ and therefore it suffices to show that $f$ is 1-1. By Proposition \ref{mbeer}, we have $f(p)=p\cap \bL$. Assume $p\ne q$ and fix $a,b$ such that $a\in p$, $b\in q$ and $a\cap b=\emptyset$. By Lemma \ref{wbeer}, there are $a',b'\in \bL$ such that $a\subs a'$, $b\subs b'$ and $a'\cap b'=\emptyset$. Then $a'\in p\cap\bL$ and $b'\in q\cap\bL$ which shows $f(p)\ne f(q)$.
\end{pf}

\begin{lm} Let $X$ be a compact Hausdorff space and let $\bK\subs\closed X$ be a basic lattice. Assume further that $\map fYX$ is a continuous map from a compact space $Y$ and let $\map h \bK{\closed Y}$ be the induced homomorphism. Then $h$ is separative. \end{lm}

\begin{pf} Assume $b\cap \inv fa=\emptyset$. By compactness, $\img fb$ is a closed set. Since $\img fb \cap a=\emptyset$, we can find $c\in\bK$ such that $c\cap a=\emptyset$ and $\img fb\subs c$. Then $b\subs\inv fc$ which shows that $h=f^{-1}$ is separative. 
\end{pf}

We finish this section with a simple application.

\begin{tw}[Alexandrov] Every compact space is a continuous image of a compact\/ $0$-dimensional space of the same weight. \end{tw}

\begin{pf}
Let $X$ be an infinite compact space and choose a basic lattice $\bL\subs\closed X$ such that $|\bL|=\w(X)$. Let $\Be$ be the Boolean algebra generated by $\bL$. Then $|\Be|=|\bL|$ and therefore $\ult\Be$ has the same weight as $X$. Since $\bL\subs\Be$, Theorem \ref{monoepik}(a) gives a surjective continuous map $\map f{\ult\Be}{X}$.
\end{pf}

\section{Extensions of compact spaces}\label{rozszerzenia}

Let $X$ be a compact space embedded into a Tikhonov cube $[0,1]^\kappa$ and assume that we are working in a countable transitive model $M$ of (a large enough part of) ZFC.
Let $\poset$ be a forcing notion in $M$ and let $G$ be a $\poset$-generic filter over $M$.
It is natural to ask how to interpret the space $X$ in $M[G]$. The obvious idea is to take $X$ as the same set and generate the topology by using open sets from the ground model $M$.
However, this new space will usually be no longer compact. Another idea is to take the closure of $X$ in $([0,1]^\kappa)^{M[G]}$, since $X$ still consists of functions from $\kappa$ to $[0,1]$ (even though $[0,1]^{M[G]}$ may be a  proper superset of $[0,1]^M$).
Now, the closure of $X$ may depend on the embedding of $X$ into a Tikhonov cube.
In fact, this is not the case, as we show below. Knowing that a compact space can be reconstructed from any of its basic lattices of closed sets, a natural option for interpretation of $X$ in $M[G]$ is to consider a basic lattice $\bL\subs\closed X$ defined in $M$ and to consider $\ult\bL$ in $M[G]$.
We shall show below that this does not depend on the choice of $\bL$ and that this is indeed a ``natural generic extension" of $X$ with respect to $G$.
In fact, we can consider extensions of compact spaces from a ZFC model $M$ whenever $M$ is a submodel of another ZFC model $N$ (not necessarily a forcing extension of $M$).

\begin{tw}
Assume $M\subs N$ are models of (a large enough part of) $ZFC$ and $X$ is a compact space in $M$ with a basic lattice $\bL$. In $N$, the space $\ult\bL$ is homeomorphic to $\ult{\closedn MX}$, where $\closedn MX$ is the lattice of closed subsets of $X$ defined in $M$.
\end{tw}

\begin{pf} Let $\bK=\closedn MX$ and let $\map h\bL{\bK}$ be the inclusion map. Then $h$ is a separative lattice homomorphism and this property is absolute. Let $f=\ultn{N}h$. Then $f$ is a continuous surjection (Theorem \ref{monoepik}). It suffices to show that $f$ is one-to-one. Fix $p\ne q$ in $\ult\bK$. Find $a\in p$ and $b\in q$ such that $a\cap b=\emptyset$. We work in $M$: Since $\bL$ is a basic lattice for $X$, there exist $a',b'\in\bL$ such that $a\cap b'=\emptyset=a'\cap b$ and $a'\cup b'=X$. 

Now, in $N$ the elements $a,b,a',b'$ are in the same relations as they were in $M$ and $f(p)$ is the unique extension of the filter $p$; therefore $b'\notin f(p)$ and $a'\in f(p)$. By the same reason, $a'\in f(q)$ and $b'\notin f(q)$. Thus $f(p)\ne f(q)$.
\end{pf}

Now we can speak about the interpretation of a compact space $X$ in any ZFC model $M$ which contains the space $X$. More precisely, assume $X$ is a compact space in a ZFC model $M$ and let $N$ be an extension of $M$. Then $\bL=\closedn MX$ is a normal separative lattice and we can define $X^N=\ultn N\bL$, where $\ultn N\bL$ denotes the space of ultrafilters of $\bL$ defined in $N$. 
If $N$ is a generic extension of $M$, i.e. $N=M[G]$ for a $\poset$-generic filter $G$ over $M$, then we shall write $X[G]$ instead of $X^{M[G]}$ and we shall say that $X[G]$ is the {\em $G$-extension of} $X$.

Now let $\Kay$ be a class of compact spaces. Formally, $\Kay=\setof{X}{\phi(X)}$, where $\phi$ is a formula of the language of set theory such that $ZFC\proves \phi(x)\implies$ ``$x$ is a compact space". We can say that $\Kay$ is {\em absolute}
if for every ZFC models $M\subs N$, we have $N\models \phi(X^N)$, whenever $X\in \Kay^M$, i.e. whenever $M\models \phi(X)$. Actually, this is the definition of {\em upward absoluteness}. There are no interesting classes of compact spaces which are downward absolute. Indeed, let $\Kay$ be a class of spaces defined by a formula $\phi$ and assume $ZFC\proves ``X\text{ is compact metric}\implies\phi(X)"$ and $ZFC\proves $``$\Kay$ is not the class of all compact spaces". Fix a countable transitive ZFC model $M$ and fix $X\in M$ such that $M\models\neg\phi(X)$. Let $\poset$ be the forcing notion collapsing $\kappa=\w(X)$ to $\aleph_0$ and let $G$ be a $\poset$-generic filter over $M$. Then the space $X[G]$ is metrizable and therefore $M[G]\models\phi(X[G])$. 

Many classes of compact spaces are absolute for obvious reasons. On the other hand, there are classes which are not absolute and classes which are absolute by highly nontrivial reasons.

\section{Absoluteness of some subclasses of Corson compacta}\label{Corson}

In this section we prove the stability of the classes of Eberlein and Gul'ko compacta. The Corson compact space of Todor\v evi\'c \cite{To2} (see also \cite[p. 287]{To1}) is an example of a non-stable Corson compact, as it is shown in \cite{KL}. We say that $X$ is a {\em stable Corson compact} if for every forcing notion $\poset$, for every $\poset$-generic filter $G$, the extension $X[G]$ is Corson.

\begin{problem}
Describe the class of stable Corson compacta without using the language of models.
\end{problem}

Recall that a space is {\em Eberlein compact} if it is homeomorphic to a compact subset of a Banach space endowed with the weak topology. We use the ``covering" characterization of Eberlein compacta, due to Rosenthal \cite{Ro}: a compact space $X$ is Eberlein if and only if there exists a family $\Yu$ consisting of open $F_\sig$ subsets of $X$ which is $T_0$ separating (i.e. $x\ne y$ implies $|\dn xy\cap U|=1$ for some $U\in\Yu$) and such that $\Yu=\bigcup_{\ntr}\Yu_n$, where each $\Yu_n$ is point-finite.
Note that if $X$ is a $0$-dimensional Eberlein compact then we may assume that the family $\Yu$ from the definition consists of clopen sets. Indeed, replace each $U\in\Yu$ by a sequence of clopen sets $\sett{W_n(U)}{\ntr}$ such that $\bigcup_{\ntr}W_n=U$ and define $\Wu_{n,k}=\setof{W_k(U)}{U\in\Yu}$. Then $\Wu=\bigcup_{n,k<\omega}\Wu_{n,k}$ is $T_0$ separating and every $\Wu_{n,k}$ is point finite.

We need the following well known fact; for the sake of completeness we give a lattice-theoretic proof.

\begin{lm}\label{zero1}
Every Eberlein compact space is a continuous image of a $0$-dimensional Eberlein compact of the same weight.
\end{lm}

\begin{pf}
Let $X$ be Eberlein compact and let $\Yu=\bigcup_{\ntr}\Yu_n$ be a $T_0$ separating collection of open $F_\sig$ subsets of $X$ such that each $\Yu_n$ is point-finite.
For each $U\in\Yu$ choose an increasing sequence of closed sets $F_n(U)$ such that $U = \bigcup_{\ntr}F_n(U)$ and $F_n(U)\subs\Int F_{n+1}(U)$ for each $\ntr$. Let $\bL$ be the lattice generated by 
$$\Ef=\setof{F_n(U)}{\ntr,\;U\in\Yu}\cup\setof{X\setminus\Int F_n(U)}{\ntr,\;U\in\Yu}.$$
Observe that $\bL$ is a basic lattice, since if $x\ne y$ in $X$ and $U\in\Yu$ is such that $x\in U$, $y\notin U$ then $x\in \Int F_n(U)$ for some $\ntr$ and consequently $A:=F_n(U)$ and $B:=X\setminus \Int F_n(U)$ are elements of $\Ef$ such that $x\notin B$, $y\notin A$ and $A\cup B=X$.

Let $\Be$ be the Boolean algebra generated by $\Ef$. Then $\bL\subs \Be$ and therefore $\ult\Be$ maps onto $X=\ult\bL$.
Refining $\Ef$, we may assume that $|\Ef|=\w(X)$, obtaining that $\ult\Be$ has the same weight as $X$. 
To see that $\ult\Be$ is Eberlein, define $\Wu=\bigcup_{k,n<\omega,\; i<2}\Wu_{k,n,i}$, where $\Wu_{k,n,0}=\setof{F_k(U)}{U\in\Yu_n}$ and $\Wu_{k,n,1}=\setof{\Int F_k(U)}{U\in\Yu_n}$ and observe that $\Wu$ generates $\Be$, i.e. it is $T_0$ separating on $\ult\Be$. Finally, each $\Wu_{k,n,i}$ is point-finite, since it contains no infinite centered subfamilies.
\end{pf}

\begin{tw} The class of Eberlein compact spaces is absolute. \end{tw}

\begin{pf}
Fix an Eberlein compact space $X$. Let $\map fYX$ be a continuous surjection, where $Y$ is a $0$-dimensional Eberlein compact space.
Let $\Yu=\bigcup_{\ntr}\Yu_n$ be a $T_0$ separating family of clopen subsets of $Y$ such that each $\Yu_n$ is point-finite. Let $\Be$ be the algebra of clopen subsets of $Y$. Being a continuous image of $\ult\Be$ is an absolute property, so it suffices to show that being a $0$-dimensional Eberlein compact is absolute. For this aim, note that in case of a family of clopen sets, ``being $T_0$ separating" is equivalent to ``generating $\Be$". Next, for each $\ntr$ define $\poset(\Yu_n)$ to be the set of all centered subcollections of $\Yu_n$. Then $\Yu_n$ is point-finite if and only if $\pair{\poset(\Yu_n)}{\sups}$ is well-founded. Now, in any extension of the universe, $\pair{\poset(\Yu_n)}{\sups}$ is the same object, since being centered is an algebraic property of finite character. Hence $\Yu$ witnesses that $\ult\Be$ is absolutely Eberlein.
\end{pf}

We now turn to the class of Gul'ko compacta. Instead of invoking the original definition, we shall recall a ``covering" characterization of Gul'ko compacta, due to Sokolov \cite{So}. For this aim, we need a definition. A family $\Yu$ of subsets of a space $X$ is {\em weakly $\sig$-point-finite} if $\Yu$ can be written as $\bigcup_{\ntr}\Yu_n$ so that for each $x\in X$, $\Yu=\setof{\Yu_n}{\ord(x,\Yu_n)<\ctbl}$, where $\ord(x,\Vee)=|\setof{V\in\Vee}{x\in V}|$.
A compact space $X$ is {\em Gul'ko compact} if and only if there exists a $T_0$ separating family of open $F_\sig$ subsets of $X$ which is weakly $\sig$-point-finite, see \cite{So}.

\begin{lm}\label{wiec} Let $X$ be a $0$-dimensional compact space. Then $X$ is Gul'ko compact iff there exists a $T_0$ separating family $\Yu$ consisting of clopen subsets of $X$ which is weakly $\sig$-point-finite. \end{lm}

\begin{pf}
Let $\Yu=\bigcup_{\ntr}\Yu_n$ be a weakly $\sig$-point-finite family of open $F_\sig$ sets which is $T_0$ separating. For each $U\in\Yu$ choose a sequence of clopen sets $W_0(U)\subs W_1(U)\subs\dots$ whose union is $U$. Let $\Wu=\bigcup_{n,k<\omega}\Wu_{k,n}$, where $\Wu_{k,n}=\setof{W_k(U)}{U\in\Yu_n}$. 
Then $\Wu$ is $T_0$ separating and consists of clopen sets.

Fix $x\in X$ and $W\in\Wu$. Then $W=W_k(U)$, where $U\in\Yu_n$ and there is $m<\omega$ such that $U\in\Yu_m$ and $\ord(x,\Yu_m)<\ctbl$. Thus also $\ord(x,\Wu_{k,m})<\ctbl$ and $W\in \Wu_{k,m}$. Taking any bijection $\map h\omega{\omega\times\omega}$ and writing $\Wu_n=\Wu_{h(n)}$, we see that $\Wu$ is weakly $\sig$-point-finite.
\end{pf}

We are going to describe a property of a Boolean algebra which is equivalent to the fact that its Stone space is Gul'ko compact. For this aim we need to define some new objects. Fix a Boolean algebra $\Be$, fix $\Yu\subs\Be$ and fix a function $\map\phi{\Yu}{\power\omega\setminus\sn\emptyset}$. Given $u\in \Yu$, define 
$$\poset(\Yu,\phi,u)=\setof{\pair sk \in\fin \Yu\times\omega}{s\text{ is centered}}.$$
Next, define a strict partial order $<$ on $\poset(\Yu,\phi,u)$ by
$$\pair sk < \pair tl\equiv s\subs t\Land k<l\Land (\forall\;m\in\phi(u)\cap k)(\exists\; v\in t\setminus s\/)\; \phi(v)=m.$$
(Recall that $k=\{0,1,\dots,k-1\}$.)

\begin{lm}\label{polityka} Let $\Be$ be a Boolean algebra. Then $\ult\Be$ is Gul'ko compact if and only if $\Be$ is generated by a family $\Yu$ such that for some function $\map\phi{\Yu}{\power\omega\setminus\sn\emptyset}$, for each $u\in \Yu$ the poset $\pair{\poset(\Yu,\phi,u)}{>}$ is well-founded.
\end{lm}

\begin{pf}
Assume first that $\ult\Be$ is Gul'ko compact. By Lemma \ref{wiec}, there is $\Yu\subs\Be$ which is $T_0$ separating and weakly $\sig$-point-finite (we identify elements of $\Be$ with the corresponding clopen sets in $\ult\Be$). Let $\Yu=\bigcup_{\ntr}\Yu_n$, where the enumeration witnesses that $\Yu$ is weakly $\sig$-point-finite. Let $\phi(u)=\setof{\ntr}{u\in\Yu_n}$. Fix $u\in\Yu$. We claim that $\poset(\Yu,\phi,u)$ is well-founded with respect to $>$. Suppose otherwise and let $\sett{\pair{s_n}{k_n}}{\ntr}$ be a $<$-increasing sequence. Then $s_\infty=\bigcup_{\ntr}s_n$ is centered so there is $p\in\ult\Be$ with $p\in\bigcap s_n$ for each $\ntr$. Find $m<\omega$ such that $u\in\Yu_m$ and $\ord(p,\Yu_m)<\ctbl$. 
Let $n_0$ be such that $k_{n_0}>m$. By the definition of $<$, for each $n\goe n_0$ there is $v_n\in (s_{n+1}\setminus s_n)\cap\Yu_m$. Then $p\in\bigcap_{n\goe n_0}v_n$ and $\setof{v_n}{\ntr}\subs\Yu_m$, so $\ord(p,\Yu_m)\goe\ctbl$, a contradiction.

Assume now that $\Yu\subs \Be$ generates $\Be$ and there is $\phi$ such that for each $u\in\Yu$, $\pair{\poset(\Yu,\phi,u)}{>}$ is well-founded.
Let $\Yu_m=\setof{u\in\Yu}{m\in\phi(u)}$. We claim that $\Yu$ is weakly $\sig$-point-finite and this fact is witnessed by $\Yu=\bigcup_{\ntr}\Yu_n$.
Suppose this is not the case and find $x\in \ult\Be$ and $u\in\Yu$ such that $\ord(x,\Yu_m)\goe\ctbl$ whenever $m\in\phi(u)$. We shall find a $<$-increasing sequence $\sett{\pair{s_n}{k_n}}{\ntr}$ in $\poset(\Yu,\phi,u)$, deriving a contradiction. 

Define $s_0=\emptyset$, $k_0=0$. Assume $\pair{s_n}{k_n}$ has been defined so that $x\in\bigcap s_n$. Fix $k>k_n$. For each $m\in\phi(u)\cap k$ find $v_m\in\Yu_m\setminus s_n$ such that $x\in v_m$. This is possible, since for $m\in\phi(u)$, $\ord(x,\Yu_m)$ is infinite and $s_n$ is finite. Define
$s_{n+1}=s_n\cup\setof{v_m}{m\in\phi(u)\cap k}$ and $k_{n+1}=k$. Then $\pair{s_n}{k_n}<\pair{s_{n+1}}{k_{n+1}}$ and $x\in \bigcap s_{n+1}$. 
This finishes the proof. 
\end{pf}

\begin{lm}\label{zero2} Every Gul'ko compact space is a continuous image of a $0$-dimensional Gul'ko compact of the same weight. \end{lm}

\begin{pf} Fix a Gul'ko compact space $X$ and let $\Yu$ witness this fact.
Define $F_n(U)$, $\Ef$, $\bL$ and $\Be$ as in the proof of Lemma \ref{zero1}. It remains to show that $\ult\Be$ is Gul'ko compact. Define $\Wu=\bigcup_{k,n<\omega,\;i<2}\Wu_{k,n,i}$ as in the proof of Lemma \ref{zero1}. Then $\Be$ is generated by $\Wu$. We are going to use Lemma \ref{polityka}. Fix a bijection $\map h\omega{\omega\times\omega\times 2}$ and define $\phi(w)=\setof{m<\omega}{w\in\Wu_{h(m)}}$.
Fix $w\in \Wu$. We shall show that
$$\pair{\poset(\Wu,\phi,w)}{>}$$
is well-founded, which will complete the proof.

Suppose $\sett{\pair{s_n}{k_n}}{\ntr}$ is a $<$-increasing sequence in $\poset(\Wu,\phi,w)$ and let $s_\infty=\bigcup_{\ntr}s_n$. By compactness, there is $x\in X$ such that $x\in\cl v$ for every $v\in s_\infty$. By the definition of $\Wu$, there are $n,k<\omega$ and $U\in\Yu_n$ such that either $w=F_k(U)$ or $w=\Int F_k(U)$. Set $i=0$ if $w=F_k(U)$ and $i=1$ otherwise.
Since $\Yu$ is weakly $\sig$-point-finite, we may assume that $\ord(x,\Yu_n)$ is finite. Let $m<\omega$ be such that $h(m)=\triple kni$. Then $m\in\phi(w)$.

Let $n_0$ be such that $k_{n_0}>m$. For each $l\goe n_0$ pick $G_l\in (s_{l+1}\setminus s_l)\cap\Wu_{k,n,i}$ (which exists by the definition of $<$). Then either $G_l=F_k(U_l)$ or $G_l=\Int F_k(U_l)$ for some $U_l\in\Yu_n$. Observe that the set $\setof{U_l}{l\goe n_0}$ must be infinite, because the sequence $\setof{G_l}{l\goe n_0}$ is one-to-one. Furthermore $x\in \bigcap_{l\goe n_0}U_l$, which shows that $\ord(x,\Yu_n)$ is infinite, a contradiction.
\end{pf}

We are now ready to prove the announced absoluteness result.

\begin{tw} The class of Gul'ko compact spaces is absolute. \end{tw}

\begin{pf} Let $X$ be Gul'ko compact. Using Lemma \ref{zero2}, find a $0$-dimensional Gul'ko compact $Y$ and a continuous surjection $\map fYX$. Let $\Be$ be the clopen algebra of $Y$. Being a continuous image of $\ult\Be$ is absolute, so it remains to observe that $\ult\Be$ is Gul'ko compact in any extension of the universe, since Lemma \ref{polityka} gives an absolute condition for this property, namely the generating set $\Yu$ and the function $\phi$ such that $\pair{\poset(\Yu,\phi,u)}{>}$ is well-founded for every $u\in\Yu$ ---all these objects remain the same in any extension of the universe and ``being a well-founded relation" is absolute. 
\end{pf}


\begin{thebibliography}{99}

%\bibitem{AP}  {\sc Alster, K.; Pol, R.}, {\em On function spaces of compact subspaces of $\Sigma$-products of the real line\/}, Fund. Math. 107 (1980), no. 2, 135--143.

\bibitem{Arhang1969} {\sc Arhangel'ski\u\i, A. V.}, 
{\em An approximation of the theory of dyadic bicompacta\/}, (Russian) Dokl. Akad. Nauk SSSR 184 (1969) 767--770.

\bibitem{Bandlow} {\sc I. Bandlow}, \textit{On the absoluteness of
openly-generated and Dugundji spaces\/}, Acta Univ. Carolin. Math. Phys. {33}
(1992), no. 2, 15--26.

\bibitem{Birkhoff} {\sc Birkhoff, G.}, {\em Lattice theory\/}. Third edition. American Mathematical Society Colloquium Publications, Vol. XXV American Mathematical Society, Providence, R.I. {\bf 1967}.

%\bibitem{BKT} {\sc Burke, M.; Kubi\'s, W.; Todor\v cevi\'c, S.}, {\em Kadec norms on spaces of continuous functions\/}, preprint.

\bibitem{Engelking} {\sc Engelking, R.}, {\em General Topology\/}. Translated from the Polish by the author. Second edition. Sigma Series in Pure Mathematics, 6. Heldermann Verlag, Berlin, {\bf 1989}.

\bibitem{Gratzer} {\sc Gr\"atzer, G.}, {\em General Lattice Theory\/}. Second edition. New appendices by the author with B. A. Davey, R. Freese, B. Ganter, M. Greferath, P. Jipsen, H. A. Priestley, H. Rose, E. T. Schmidt, S. E. Schmidt, F. Wehrung and R. Wille. Birkhauser Verlag, Basel, {\bf 1998}.

%\bibitem{Kelley} {\sc Kelley, J.}, {\em General Topology\/}. Reprint of the 1955 edition [Van Nostrand, Toronto, Ont.]. Graduate Texts in Mathematics, No. 27. Springer-Verlag, New York-Berlin, {\bf 1975}.

%\bibitem{K2004a} {\sc Kubi\'s, W.}, {\em Compact spaces with many retractions\/}, preprint.

\bibitem{KL} {\sc W. Kubi\'s, A. Leiderman}, {\em Semi-Eberlein spaces\/}, Topology Proc. 28 (2004), No. 2, 603--616 

\bibitem{KOSz} {\sc W. Kubi\'s, O. Okunev, P. Szeptycki}, {\em On some classes of Lindelof Sigma-spaces}, Topology Appl. Vol. 153 (2006) 2574--2590.

%\bibitem{vanMill1982} {\sc van Mill, J.}, {\em A homogeneous Eberlein compact space which is not metrizable\/}, Pacific J. Math. 101 (1982), no. 1, 141--146.

\bibitem{Ro} {\sc H. Rosenthal}, {\em The heredity problem for weakly compactly generated Banach spaces\/}, Compositio Math. 28 (1974), 83--111.

%\bibitem{Rudin1956} {\sc Rudin, W.}, {\em Homogeneity problems in the theory of \v Cech compactifications\/}, Duke Math. J. 23 (1956), 409--419.

\bibitem{So} {\sc G. A. Sokolov}, {\em On some classes of compact spaces lying in $\Sigma$-products\/}, Comment. Math. Univ. Carolin. 25 (1984), no. 2, 219--231.

%\bibitem{Tkachuk1994} {\sc Tkachuk, V.}, {\em A glance at compact spaces which map ``nicely" onto the metrizable ones\/}, Topology Proc. 19 (1994), 321--334.

\bibitem{To2} {\sc S. Todor\v cevi\'c}, {\em Stationary sets, trees and continuums\/}, Publ. Inst. Math. (Beograd) \textbf{27} (1981), pp. 249--262.

\bibitem{To1} {\sc S. Todor\v cevi\'c}, {\em Trees and linearly ordered sets\/}, Handbook of set-theoretic topology, 235--293, North-Holland, Amsterdam, 1984.

\bibitem{TodorcevicBaireClassOne} {\sc S. Todor\v cevi\'c}, {\it Compact
subsets of the first Baire class\/}, J. Amer. Math. Soc. \textbf{12} (1999),
no. 4, 1179--1212.

\end{thebibliography}
\end{document}